\documentclass[11pt]{amsart}
\usepackage{amsmath,amssymb,amsfonts,graphics}

\newtheorem{thm}{Theorem}
\newtheorem{lem}[thm]{Lemma}

\newtheorem{prop}[thm]{Proposition}

\theoremstyle{definition}
\newtheorem{defn}[thm]{Definition}

\newtheorem{rmk}[thm]{Remark}

\newcommand{\R}{\mathbb{R}}
\newcommand{\N}{\mathbb{N}}

\newcommand{\D}{\mathcal{D}}
\newcommand{\fA}{\mathcal{A}}

\newcommand{\fT}{\mathcal{T}}
\newcommand{\Del}{\Delta}
\newcommand{\Inn}{{\rm Inn}}
\newcommand{\supp}{{\rm supp}}
\newcommand{\spn}{{\rm span}}
\newcommand{\cross}{\times}

\begin{document}

\title[Smooth and weak synthesis of the anti-diagonal]
{Smooth and weak synthesis of the anti-diagonal in Fourier algebras of Lie groups}

\author{B. Doug Park}
\address{Department of Pure Mathematics, University of Waterloo, Waterloo,
Ontario, N2L 3G1, Canada}
\email{bdpark@math.uwaterloo.ca}

\author{Ebrahim Samei}
\address{Department of Mathematics and Statistics, University of Saskatchewan, Saskatoon, Saskatchewan, S7N 5E6, Canada}
\email{samei@math.usask.ca}




\begin{abstract}
Let $G$\/ be a Lie group of dimension $n$, and let $A(G)$ be the
Fourier algebra of $G$. We show that the anti-diagonal
$\check{\Delta}_G=\{(g,g^{-1})\in G\times G \mid g\in G\}$ is
both a set of local smooth synthesis and a set of local weak
synthesis of degree at most $[\frac{n}{2}]+1$ for $A(G\times G)$.
We achieve this by using the concept of the cone property in
\cite{ludwig-turowska}. For compact $G$, we give an alternative
approach to demonstrate the preceding results by applying the
ideas developed in \cite{forrest-samei-spronk}.  We also present
similar results for sets of the form $HK$, where both $H$\/ and $K$\/
are subgroups of $G\times G\times G\times G$\/ of diagonal forms.
Our results very much depend on both the geometric and the
algebraic structure of these sets.
\end{abstract}

\maketitle

\section{Introduction}

The concept of weak (spectral) synthesis was introduced and
studied by Warner in \cite{W}. This is a generalization of the
notion of sets of synthesis for a regular Banach algebra of
continuous functions. Motivation arose in \cite{W} from studying
the question of when the union of two sets of synthesis is a set
of synthesis. This concept also appeared in the earlier work of
others such as Varopoulos \cite{V} who re-proved the well-known
result of L. Schwartz: the $(n-1)$-dimensional sphere $S^{n-1}$
is a set of weak synthesis for the Fourier algebra $A(\R^n)$ of
degree $[\frac{n}{2}+1]$ (see also \cite{Kan} and the reference
therein for more examples).

There is another notion of synthesis, known as smooth synthesis,
which is considered for regular Banach algebras of continuous
functions on $C^\infty$ manifolds \cite{ludwig-turowska} (see
also \cite{guo} and \cite{KM}). This property is, roughly
speaking, a suitable concept when the set is a smooth curve. We
refer the reader to Domar's survey article \cite{D3} for more
information. An important observation by Domar \cite{D2} showed
that smoothness of the curve alone is not sufficient to imply its
smooth synthesis. However, Domar showed that if a desirable curve
in $\R^n$ also satisfies certain curvature condition, then both smooth
synthesis and weak synthesis hold \cite{domar}. This curvature
condition was generalized and modified in \cite{ludwig-turowska}
(see also \cite{guo} and \cite{KM}). A closed set satisfying this
condition is said to have {\it the cone property} (see Definition
\ref{defn: cone}).

Let $G$\/ be a locally compact group.  In recent years, owing much
of it to the applications of operator space theory to harmonic analysis,
the structure of the anti-diagonal
$\check{\Delta}_G=\{(g,g^{-1})\in G\times G \mid g\in G\}$ became
closely related to the cohomological property of Fourier algebra
$A(G)$. For instance, it is shown in \cite{FR} that $A(G)$ is
amenable exactly when $\check{\Delta}_G$ is an element of the
coset ring of $G\cross G$, or exactly when $G$\/ admits an abelian
subgroup of finite index. In an analogous result, it is given in
\cite[Theorem 2.4]{FSS2} for a compact group $G$, a full
characterization of when $A(G)$ is weakly amenable: when the
connected component of the identity $G_e$ is abelian, or
equivalently, $\check{\Delta}_G$ is a set of synthesis for
$A(G\times G)$ (see also \cite[Theorem~3.7]{forrest-samei-spronk}). This can be compared to the
well-known result that closed subgroups are of synthesis for
Fourier algebras \cite[Theorem 3]{TT} and $\check{\Delta}_G$ is a
subgroup of $G\times G$\/ exactly when $G$\/ is abelian. Hence the
preceding result shows that the synthesis property of
$\check{\Delta}_G$ is inherited from the group structure of $G_e$
or the cohomological behavior of $A(G)$.

In this article, for a Lie group $G$, we investigate both smooth
and weak synthesis of certain diagonal-type sets including the
anti-diagonal. The paper is organized as follows.

In Section \ref{S:prim}, we state the general background of the
concepts and notions required. In Section \ref{S:General} we
summarize and slightly modify the main results in \cite[Section
4]{ludwig-turowska}. This, in particular, is essential for the
results in Section \ref{Diag type-smooth-weak} with regard to \emph
{local}\/ smooth and weak synthesis.

In Subsection \ref{S:anti-diag}, we first show that
$\check{\Delta}_G$ is a smooth submanifold of $G\times G$\/ with
the cone property. Then by using the tools developed in
\cite{ludwig-turowska} and their modification in Section
\ref{S:General}, we show that $\check{\Delta}_G$ is both a set of
local smooth synthesis and a set of local weak synthesis for
$A(G\times G)$. Moreover the degree of the nilpotency is at most
$[\frac{n}{2}]+1$ where $n$\/ is the dimension of $G$ (Theorem
\ref{T:smooth weak syn-anti diagonal}). The word ``local" is
redundant if $A(G)$ has an approximate identity. However since it
is not known whether this holds in general, we do not know
whether different notations of synthesis and their corresponding
local synthesis coincide.

It is shown in \cite[Corollary 3.2]{forrest-samei-spronk} that if
$G$\/ is a compact, connected, non-abelian Lie group, then
$(\Delta_G\times \Delta_G)\Delta_{G\times G}$ is a set of
non-synthesis for $A(G\times G\times G\times G)$. This gives us an
example of a set of non-synthesis which is the product of two closed
subgroups. In Subsection \ref{S:product diag}, we show that, for a
Lie group $G$, $(\Delta_G\times \Delta_G)\Delta_{G\times G}$ is both
a set of local smooth and local weak synthesis for $A(G\times
G\times G\times G)$. Moreover the degree of the nilpotency is at
most $[\frac{3n}{2}]+1$ where $n$ is the dimension of $G$ (Theorem
\ref{T:smooth weak syn-product}). This is done in a similar fashion
to that of Subsection \ref{S:anti-diag}, i.e. by showing that
$(\Delta_G\times \Delta_G)\Delta_{G\times G}$ is a smooth
$3n$-dimensional submanifold of $G\times G\times G\times G$\/ with
the cone property.

The rest of the paper is devoted to presenting an alternative
approach to prove the preceding results in the case of a compact
Lie group $G$. In Section \ref{S:Proj-smooth-weak}, we prove the
projection theorem for smooth and weak synthesis for an arbitrary
closed subgroup $K$ (Theorems \ref{T:coset-weak syn} and
\ref{T:coset-smooth syn}). The main theory required to achieve
this is developed in \cite[Sections 1.3 and 2]
{forrest-samei-spronk}. We apply these results in Section
\ref{S:Alternative} to obtain the smooth and weak synthesis of
$\check{\Delta}_G$ and $(\Delta_G\times \Delta_G)\Delta_{G\times
G}$ (Theorem \ref{T:smooth weak syn-anti
diagonal-product-alternative}).

\section{Preliminaries and notations}
\label{S:prim}

\subsection{Notations of Synthesis}
Let $\fA$\/ be a Banach algebra contained in $C_0(X)$ for some
locally compact Hausdorff space $X$. We define for any closed
subset $E$\/ of $X$\/
\begin{align*}
I_\fA(E)&=\{f\in\fA \mid f(x)=0 \ \text{for all } x\in E\},  \\
I_\fA^0(E)&=\{f\in\fA \mid \supp{f}\cap E=\varnothing,\
\supp{f}\text{ is compact}\},\\
J_\fA(E)&=\overline{\{a \in I_\fA(E) \mid \supp\, a \ \text{is compact}
\}},
\end{align*}
where $\supp{f}= \overline{\{x\in X \mid f(x)\not=0\}}$ and the closure in
$J_\fA(E)$ is taken with respect to the norm $\|\cdot \|_\fA$. If $X$ is, in
addition, a smooth manifold, then we let $\D(X)$ be the space of
all compactly supported $C^\infty$ functions on $X$ and denote by
$J^\D_\fA(E)$ the closure (in $\fA$) of the space of all elements in
$\D(X)\cap \fA$ which vanishes on $E$. When there is no fear of
ambiguity, we write $I(E)$ instead of $I_\fA(E)$, $I_0(E)$
instead of $I^0_\fA(E)$, $J(E)$ instead of $J_\fA(E)$, and
$J_\D(E)$ instead of $J^\D_\fA(E)$.

Suppose that $\fA$\/ is regular on $X$. We say that $E$\/ is {\it a
set of synthesis} ({\it local synthesis}) for $\fA$\/ if $I_0(E)$ is
dense in $I(E)$ ($J(E)$). More generally, we say that $E$\/ is a
set of \emph{weak synthesis $($local weak synthesis$)$ for $\fA$
of degree at most $d$}\/ if $I(E)^d=\overline{I_0(E)}$
($J(E)^d=\overline{I_0(E)}$) for some positive integer $d$, i.e.
if $I(E)/\overline{I_0(E)}$ ($J(E)/\overline{I_0(E)}$) is
nilpotent of degree at most $d$. If $X$\/ is, in addition, a smooth
manifold, then $E$\/ is {\it a set of smooth synthesis}\/ for $\fA$\/ if
$J_\D(E)=I(E)$ and it is {\it a set of local smooth synthesis}\/
for $\fA$\/ if $J_\D(E)=J(E)$ \cite{ludwig-turowska} (see also
\cite{guo} and \cite{KM}). It is clear that every set of (weak)
synthesis is a set of local (weak) synthesis and the converse
holds if $\fA$\/ has an approximate identity with compact support.

Suppose further that $X$\/ is the Gelfand spectrum of $\fA$. Then
$I(E)$ is the largest and $I_0(E)$ is the smallest ideal in $\fA$
whose hull is $E$ \cite[Proposition 4.1.20]{D}. Thus $E$\/ is a set
of synthesis for $\fA$\/ if and only if there is a unique closed
ideal in $\fA$\/ whose hull is $E$. Let $\fA_c$ be the set of all
elements in $\fA$\/ with compact support. If $\fA_c$ is dense in
$\fA$, i.e. $\fA$\/ is a \emph{Tauberian algebra} \cite{Ric}, then
$J(E)$ is the maximal ideal of $\fA$\/ having $E$\/ as its hull and
being essential as a Banach $\fA$-bimodule. So if $E$\/ is a set of
local synthesis, then $J(E)$ is the only closed ideal in $\fA$\/
with this property.

\subsection{Fourier algebras}
Let $G$\/ be a locally compact group with a fixed left Haar measure.
Given a function $f$\/ on $G$\/ the left and right translation of $f$\/
by $x\in G$\/ is denoted by $(L_xf)(y)=f(xy)$ and
$(R_xf)(y)=f(yx)$, respectively. Let $P(G)$ be the set of all
continuous positive definite functions on $G$\/ and let $B(G)$ be its
linear span. The space $B(G)$ can be identified with the dual of
the group $C^*$-algebra $C^*(G)$, this latter being the
completion of $L^1(G)$ under its largest $C^*$-norm. With
pointwise multiplication and the dual norm, $B(G)$ is a
commutative regular semisimple Banach algebra. The Fourier
algebra $A(G)$ is the closure of $B(G)\cap C_c(G)$ in $B(G)$. It
was shown in \cite{Em} that $A(G)$ is a commutative regular
semisimple Banach algebra whose carrier space is $G$. Also, if
$\lambda$ is the left regular representation of $G$\/ on $L^2(G)$
then, up to isomorphism, $A(G)$ is the unique predual of $VN(G)$,
the von Neumann algebra generated by the representation $\lambda$.

\section{general results on smooth and weak
synthesis}\label{S:General}

In \cite[Section 4]{ludwig-turowska}, the smooth and weak
synthesis properties were studied in the Fourier algebra $A(G)$ of a
Lie group $G$. In this section, we summarize and slightly modify
their main results. These modifications will be used to obtain
some of the main results of our article.

\subsection{Connection between smooth and weak synthesis}
We start with the following theorem which is a generalization of
Theorem 4.3 and Corollary 4.4 in \cite{ludwig-turowska}.

\begin{thm}\label{T:smooth-weak synthesis}
Let $G$ be a Lie group.  Let $M$ be a smooth $m$-dimensional
submanifold of $G$, and let $E\subseteq M$ be closed in
$G$. Then:\\
{\rm (i)} $J_\D(E)^{[\frac{m}{2}]+1}=\overline{I_0(E)}$.\\
{\rm (ii)} If $E$ is a set of smooth synthesis, then it is a set of weak
synthesis of degree at most $[\frac{m}{2}]+1$.\\
{\rm (iii)} If $E$ is a set of local smooth synthesis, then it is a set of local weak
synthesis of degree at most $[\frac{m}{2}]+1$.
\end{thm}

\begin{proof}
(i) For simplicity, let $d=[\frac{m}{2}]+1$. It suffices to show that $J_\D(E)^d \subseteq \overline{I_0(E)}$.
Let $f\in J_\D(E)$ with compact support. We will follow a similar argument to that
of \cite[proposition 1.6]{Kan} to demonstrate that $f^d\in \overline{I_0(E)}$. Using
the partition of unity, it is sufficient to show that $f^d$ is ``locally
in $\overline{I_0(E)}$" i.e. for every $x\in G$, there is a neighborhood
$U_x$ of $x$\/ in $G$\/ and an element $f_x \in \overline{I_0(E)}$ such that
$f^d=f_x$ on $U_x$. Since $A(G)$ is regular, $f^d$ belongs locally to $\overline{I_0(E)}$
at every point $x\in G\setminus E$. On the other hand, for every $x\in E$, let $V_x$ be a compact neighborhood
of $x$\/ in $G$, and let $E_x=E\cap V_x$. It follows from \cite[Theorem 4.3]{ludwig-turowska}
that
$$
J_\D(E_x)^d=\overline{I_0(E_x)}.
$$
In particular,
\begin{equation}\label{eq:f^d}
f^d\in \overline{I_0(E_x)}.
\end{equation}
Take $h_x\in A(G)$ such that
$\supp\hspace{1pt} h_x \subseteq V_x$ and $h_x=1$ on a neighborhood of $x$ in $G$,
say $U_x$. Let $\epsilon >0$. By (\ref{eq:f^d}), there is $g_x\in I_0(E_x)$ such that
$$\|f^d-g_x \| < \frac{\epsilon}{\|h_x\|} .$$
Thus $\|h_xf^d-h_xg_x \|< \epsilon$. Moreover, $h_xg_x=0$ on a neighborhood
of $E$. To see this, we observe that $h_x=0$ on a neighborhood of $G\setminus V_x$
and $g_x=0$ on a neighborhood of $E_x$. So $h_xg_x=0$ on a neighborhood
of $E_x\cup(G\setminus V_x)$ which contains $E$. Hence $h_xg_x\in I_0(E)$. Since
$\epsilon >0$ was arbitrary, we have
$$h_xf^d\in \overline{I_0(E)}.$$
Moreover, $h_xf^d=f^d$ on $U_x$ since $h_x=1$ on $U_x$. That is
$f$ belongs locally to $\overline{I_0(E)}$ at every point $x\in
E$. Since $\supp f$ is compact, an argument using partition of unity
implies that $f^d\in \overline{I_0(E)}$. The final result follows
from the fact that $J_\D(E)$ is the closure of $\{f \in J_\D(E)
\mid  \supp\, f \ \text{is compact} \}$.

Parts (ii) and (iii) follow immediately from part (i).
\end{proof}

We note that the preceding theorem implies that any closed subset $E$\/ of a Lie group
$G$\/ with (local) smooth synthesis has (local) weak synthesis. Moreover the degree of
the nilpotency is dominated by the dimensions of smooth submanifolds containing $E$.

\subsection{Cone property}
Motivated by Theorem \ref{T:smooth-weak synthesis}, we would like to find a sufficient condition
that implies a closed set to be of (local) smooth synthesis. One such condition is
the cone property that was introduced in \cite{domar}, \cite{guo} and \cite{KM} for the
case of $\R^n$ and in \cite{ludwig-turowska} for a general Lie group.

\begin{defn}\label{defn: affine}
(cf.\ \cite[Definition 4.5]{ludwig-turowska}) \
Given a Lie group $H$, let ${\rm Aut}(H)$ denote its (continuous) automorphism
group. It is well known that ${\rm Aut}(H)$ is also a Lie group. Also, since
every continuous homomorphism between Lie groups is analytic (see \cite[Ch.\ II, Theorem 2.6]{helgason}),
each element of ${\rm Aut}(H)$ is smooth.

{\rm (i)} \  Any pair $t\in H$ and $a\in {\rm Aut}(H)$ give rise to a
mapping $\varphi:H\rightarrow H$, $\varphi(s)= a(ts)$, which will
be called an \emph{affine transformation}\/ of $H$.

{\rm (ii)} \ We say that a Lie group $A$\/ is a \emph{group of affine transformation}
of $H$\/ if $A$\/ acts smoothly by affine transformations on $H$. Smoothly here means that
the mapping $A\times H \to H$, $(a,x) \mapsto a(x)$ is smooth.

{\rm (iii)} \ Let ${\rm Aff}(H)$ denote the group of all affine transformations of $H$.
Here we identify ${\rm Aff}(H)$ with the semidirect product $H \ltimes_\rho {\rm Aut}(H)$,
where $\rho : {\rm Aut}(H) \to {\rm Aff}(H)$ is the identity automorphism. In particular,
${\rm Aff}(H)$ is a Lie group and acts smoothly on $H$\/ by the action
${\rm Aff}(H)\times H \to H$, $((t,a),s) \mapsto a(ts)$. Moreover if $A$\/ is any
group of affine transformations of $H$, then the identity map from $A$\/ into ${\rm Aff}(H)$
is a smooth monomorphism. We note that image of $A$\/ under this map does not have to be closed in ${\rm Aff}(H)$.
\end{defn}

\begin{defn}\label{defn: cone}
(cf.\ \cite[Definition 4.6]{ludwig-turowska}) \
Let $H$\/ be a Lie group and let $M$\/ be a smooth
$m$-dimensional submanifold of $H$. We say that a subset $E$\/ of
$M$ has the \emph{cone property}\/ if the following holds:
\begin{itemize}
\item[(i)] $E$\/ is closed in $H$.

\item[(ii)] For every $x\in E$, there exists an open neighborhood $U_x$ of $x$\/ in $H$\/ and a $C^{\infty}$ mapping $\psi_x$ from an open subset $W_x\subset \R^m$ containing $0$ into a Lie group of affine transformations $A_x$ on $H$\/ such that $\psi_x(0)={\rm id}_H$ and there exists an open subset $W_x^0\subset W_x$ such that:
\begin{itemize}
\item[(a)] $0$ is in the closure of $W_x^0$.

\item[(b)] For every $y\in U_x\cap E$, $\psi_x(W_x^0)y$ is contained in $E$\/ and open in $M$ and the mapping $W_x^0\rightarrow \psi_x(W_x^0)y$, $t\mapsto \psi_x(t)y$, is a diffeomorphism.

\end{itemize}
\end{itemize}
\end{defn}

\begin{thm}\label{T:cone prop-smooth weak syn}
Let $G$ be a Lie group.  Let $M$ be a smooth $m$-dimensional submanifold of $G$,
and let $E\subseteq M$ be a subset with the cone property. Then:\\
{\rm (i)} $E$ is a set of local smooth synthesis.\\
{\rm (ii)} $E$ is a set of local weak synthesis of degree at most\/ $[\frac{m}{2}]+1$.
\end{thm}

\begin{proof}
(i) Let $T$\/ be an element of $J_\D(E)$, and let $f\in A(G)$ with compact support. Then
$f\cdot T\in VN(G)$ with compact support and vanishes $J_\D(E)$. It is shown in the
proof of \cite[Theorem 4.3]{ludwig-turowska} that $f\cdot T$ must vanish $I(E)$. Therefore
$T=0$ on $I(E)A(G)$, and so, $T$ annihilates $J(E)$. This means that $J_\D(E)=J(E)$.\\
(ii) This follows from part (i) and Theorem \ref{T:smooth-weak synthesis}(iii).
\end{proof}

One of the main consequences of the preceding theorem is
\cite[Corollary 4.9]{ludwig-turowska} where it is shown that
certain orbits are of local weak synthesis. We note that in the
statement of \cite[Corollary 4.9]{ludwig-turowska} the phrase
``weak synthesis" is used instead of ``local weak synthesis".
However a careful examination of the proof of \cite[Corollary
4.9]{ludwig-turowska} and a comparison with the argument in
Theorem \ref{T:smooth-weak synthesis} (see also \cite[Proposition
1.6]{Kan}) show that their argument actually proves the local
weak synthesis property. For that sake, we state the correct
version of that corollary below.

\begin{thm}\label{T:smooth weak syn-orbits}
Let $G$ be a Lie group, and let $B$ be a group of affine transformations of $G$.
Let $\omega\subseteq G$ be a closed $m$-dimensional $B$-orbit in $G$. Then:\\
{\rm (i)} $\omega$ is a set of local smooth synthesis.\\
{\rm (ii)} $\omega$ is a set of local weak synthesis of degree at most\/ $[\frac{m}{2}]+1$.
\end{thm}

\section{Diagonal-type sets of smooth and weak
synthesis}\label{Diag type-smooth-weak}

\subsection{Anti-diagonal}\label{S:anti-diag}

Let $G$\/ be an $n$-dimensional Lie group.  Let
$$\check{\Delta}_G=\{(g,g^{-1})\in G\times G \mid g\in G\},$$
which is called the {\it anti-diagonal}\/ of $G\times
G$. Note that $\check{\Delta}_G$ is a closed smooth
$n$-dimensional submanifold of $G\times G$\/ since it is equal to
the inverse image $\mu^{-1}(e)$, where $\mu:G\times G\rightarrow
G$\/ is the multiplication map given by $\mu(g_1,g_2)=g_1g_2$ and
the identity element $e\in G$\/ is a regular value of $\mu$.
Note that the map $\alpha: G \rightarrow \check{\Delta}_G$ given
by
\begin{equation}\label{eq: alpha}
\alpha(g)=(g,g^{-1})
\end{equation}
is a diffeomorphism whose inverse is the projection map onto the first factor.

It is shown in \cite[Theorem 2.4]{FSS2} that if $G$ is compact
with non-abelian connected component, then $\check{\Delta}_G$ is
not a set of synthesis for $A(G)$. In this section, we show that
for a general Lie group $G$, $\check{\Delta}_G$ is always a set of
local weak synthesis. In particular, it is of weak synthesis when
$A(G)$ has an approximate identity. We start with the following
lemma which shows that $\check{\Delta}_G$ has the cone property.

\begin{lem}\label{L:cone property-anti diagonal}
Let $G$ be an $n$-dimensional Lie group. Then the submanifold\/
$\check{\Delta}_G$ has the cone property in $G\times G$.
\end{lem}

\begin{proof}
We use the notation in Definition~\ref{defn: cone}.  We set
$H=G\times G$, $E=M=\check{\Delta}_G$, $m=n$, and $U_x=H=G\times G$\/
for every $x\in\check{\Delta}_G$. Let $\mathfrak{g}$ denote the Lie
algebra of $G$.  Let $W$\/ be an open subset of $\mathfrak{g}\cong
\R^n$ containing $0$ such that the exponential map restricted to
$W$, $\exp :W\rightarrow \exp(W)\subset G$, is a diffeomorphism.
We set $W_x^0=W_x=W$, and set $A_x ={\rm Aff}(G\times G)$ for
every $x\in \check{\Delta}_G$.  We write $\tilde{t}=\exp(t)$ for
short.  We define $\psi_x: W\rightarrow {\rm Aff}(G\times G)$,
also independent of $x$, as follows:  for $t\in W$ and
$(g_1,g_2)\in G\times G$,
\begin{eqnarray*}
\psi_x(t)(g_1,g_2) &=& ({\rm Inn}(\tilde{t}) \times {\rm id}_G)((\tilde{t},\tilde{t}^{-1})\cdot (g_1,g_2))\\
&=& (\tilde{t}^{-1}(\tilde{t} g_1)\tilde{t}, \tilde{t}^{-1} g_2)\\
&=& (g_1\tilde{t}, \tilde{t}^{-1}g_2),
\end{eqnarray*}
where ${\rm Inn}(\tilde{t})$ denotes the inner automorphism of
$G$\/ corresponding to $\tilde{t}\in G$. Each $\psi_x(t)$ is the
composition of a left multiplication by
$(\tilde{t},\tilde{t}^{-1})$ and an inner automorphism ${\rm
Inn}(\tilde{t}) \times {\rm id}_G \in{\rm Aut}(G\times G)$, and
hence is an affine transformation of $G\times G$. Clearly,
$\psi_x$ is smooth and $\psi_x(0)={\rm id}_{G\times G}$.  For
every $y=(g,g^{-1})\in \check{\Delta}_G$ and $t\in W$, we have
\[
\psi_x(t)y=(g\tilde{t},\tilde{t}^{-1}g^{-1})=(g\tilde{t},(g\tilde{t})^{-1})\in \check{\Delta}_G.
\]
The mapping $W\rightarrow \psi_x(W)y$ is just the composition
$\alpha \circ L_g \circ\exp$, where $\alpha$\/ is the
diffeomorphism given by (\ref{eq: alpha}). Since $\exp$ is a
diffeomorphism on $W$, the mapping $W\rightarrow \psi_x(W)y$ is
also a diffeomorphism.
\end{proof}

\begin{thm}\label{T:smooth weak syn-anti diagonal}
Let $G$ be an $n$-dimensional Lie group. Then:\\
{\rm (i)} $\check{\Delta}_G$ is a set of local smooth synthesis for $A(G)$.\\
{\rm (ii)} $\check{\Delta}_G$ is a set of local weak synthesis for $A(G)$ of degree at most\/ $[\frac{n}{2}]+1$.
\end{thm}

\begin{proof}
Our results follow from Theorem \ref{T:cone prop-smooth weak syn} and Lemma \ref{L:cone property-anti diagonal}.
\end{proof}

We note that Theorem \ref{T:smooth weak syn-anti diagonal} can
also be deduced directly from Theorem \ref{T:smooth weak
syn-orbits}. This is because $\check{\Delta}_G$ is the orbit
generated by the action $G\times (G\times G) \to G\times G$, $(x,
(g,h)) \mapsto (gx, x^{-1}h)$ and the identity element $(e,e)\in
G\times G$. As it is shown in the proof of Lemma \ref{L:cone
property-anti diagonal}, $G$\/ is a group of affine transformations
of $G\times G$ by the above action. Nevertheless, we would like
to point out that our original proof of Theorem \ref{T:smooth weak
syn-anti diagonal} is more straightforward. The proof of Theorem
\ref{T:smooth weak syn-orbits} which is \cite[Corollary
4.9]{ludwig-turowska}, although more general, is based on a
``local slicing" of an orbit to subsets with the cone property.
However we directly proved in Lemma \ref{L:cone property-anti
diagonal} that in the case of the anti-diagonal, the cone property
always holds. Moreover, as we will see in the following proposition,
any attempt toward making the anti-diagonal having the cone
property is ``locally" similar to the construction in the proof of Lemma
\ref{L:cone property-anti diagonal}.

\begin{rmk}\label{R:cone prop-connected}
In Definition \ref{defn: cone}, let $W'_x$ be the connected
component of $\{0\}$ in $W_x \subseteq \R^m$. Then $W'_x$ is an
open connected neighborhood of $0\in \R^m$. So, by replacing
$W_x$ with $W'_x$ and $W^0_x$ with $W^0_x\cap W'_x$, we see that
the assumption of Definition \ref{defn: cone} still holds. Hence,
we can assume that, for every $x\in G$, $W_x$ is connected. This,
in particular, implies that $\psi_x$ maps $W_x$ into the connected
component of the identity $(e, {\rm id}_G)$ in ${\rm Aff(G)}$.
\end{rmk}

\begin{prop}
Let $G$ be a connected, semisimple, Lie group, and let\/ $\Inn(G)$
denote the group of inner automorphisms of $G$. For every $x\in
G$, let $\psi_x$ and $W^0_x$ be the corresponding objects to
$(x,x^{-1})$ in Definition\/ $\ref{defn: cone}$ and Remark\/
$\ref{R:cone prop-connected}$ for $\check{\Delta}_G$. Then there
exist smooth maps $\rho_x : W^0_x \to G$ and $\varphi_x : W^0_x
\to \Inn(G)$ such that
$$\psi_x(r)(w,z)=(\varphi_x(r)[w\rho_x(r)^{-1}],
\varphi_x(r)[\rho_x(r)z])$$ for all\/ $r\in W^0_x$ and $w,z \in
G$. That is,
$$\psi_x(r)= [R_{\rho_x(r)^{-1}} \otimes L_{\rho_x(r)}]
\circ [\varphi_x(r) \otimes \varphi_x(r)] \quad (x\in G,\; r\in W^0_x).$$
\end{prop}

\begin{proof}
Let $E=M=\check{\Delta}_G$, $x\in G$, and $U_x=U_{(x,x^{-1})}$ be as in
Definition \ref{defn: cone} for $\check{\Delta}_G$. There is a
symmetric compact neighborhood $V$ of the identity $\{e\}$ in $G$
such that
$$(xV\times Vx^{-1}) \subseteq U_x.$$
Let $v\in V$. Then we have
\begin{equation}\label{eq:prop-1}
(xv, v^{-1}x^{-1})=(xv, (xv)^{-1})\in U_x\cap M.
\end{equation}
Let $\sigma=\psi_x(r) \in \psi_x(W^0_x)$. Since $G$\/ is semisimple
and connected, the connected component of the identity $((e,e),
{\rm id}_{G\times G})$ in ${\rm Aff}(G\times G)$ is $G\times G\times
\Inn(G\times G)$. Therefore, by Definition~\ref{defn: affine} and Remark~\ref{R:cone
prop-connected}, there are $a=b\otimes c \in \Inn(G\times G)$ and
$(t_1,t_2)\in G\times G$ such that
$$\sigma=L_{(t_1,t_2)}a.$$
Thus, from (\ref{eq:prop-1}) and Definition
\ref{defn: cone},
$$(b(t_1xv), c(t_2v^{-1}x^{-1}))=a(t_1xv, t_2v^{-1}x^{-1})=\sigma(xv,v^{-1}x^{-1}) \in
M.$$
That is,
\begin{equation*}
b(t_1xv)=c(xvt^{-1}_2), \ \ \forall v\in V.
\end{equation*}
In particular,
\begin{equation}\label{eq:prop-3}
b(t_1x)=c(xt^{-1}_2).
\end{equation}
Thus we have
\begin{eqnarray*}
c(x)c(v)c(t_2^{-1}) & = & c(xvt_2^{-1})
= b(t_1xv) \\
&=& b(t_1x)b(v)
= c(xt_2^{-1})b(v) \\
&=& c(x)c(t_2^{-1})b(v).
\end{eqnarray*}
Therefore $b=c\circ \Inn(t_2^{-1})$ on $V$, and so, $b=c\circ
\Inn(t_2^{-1})$ on the closed subgroup $H$\/ generated by $V$. Since
$V$ is open, so is $H$. Hence $H=G$\/ since $G$\/ is connected. That
is,
\begin{equation}\label{eq:prop-4}
b=c\circ \Inn(t_2^{-1}).
\end{equation}
Now, from (\ref{eq:prop-3}) and (\ref{eq:prop-4}), it follows that
\begin{equation*}
b(t_1x) = c(xt_2^{-1})
= (c\circ \Inn(t_2^{-1}))(t_2^{-1}x)
= b(t_2^{-1}x).
\end{equation*}
So $t_1x=t_2^{-1}x$ which implies that
$$t_1=t_2^{-1}.$$
Therefore, for every $s,s'\in G$,
\begin{eqnarray*}
\sigma(s,s') & = & a(t_1s, t_2s') \\
& = & ((c\circ \Inn(t_2^{-1}))(t_2^{-1}s), c(t_2s')) \\
& = & (c(st_2^{-1}), c(t_2s')).
\end{eqnarray*}
Now define $\rho_x(r):=t_2$ and $\varphi_x(r):=c$.  Then $\rho_x :
W^0_x \to G$\/ and $\varphi_x : W^0_x \to \Inn(G)$ are well-defined
and
$$\psi_x(r) = \sigma = [R_{\rho_x(r)^{-1}} \otimes L_{\rho_x(r)}]
\circ [\varphi_x(r) \otimes \varphi_x(r)].$$
Moreover, since $\psi_x$ is smooth, the above equality implies
that both $\rho_x$ and $\varphi_x$ are smooth.
\end{proof}

\subsection{Certain products of diagonals}\label{S:product diag}

Let $G$ be a Lie group, and let
$$\Delta_G=\{(g,g)\in G\times G \mid g\in G\},$$
which is the {\it diagonal}\/ of $G\times G$. It is well-known
that closed subgroups are sets of synthesis for Fourier algebras
\cite[Theorem 3]{TT} (see also \cite{H1}). It is shown in
\cite{forrest-samei-spronk} that this could fail if we consider
products of subgroups. More precisely, it is shown in
\cite[Corollary 3.2]{forrest-samei-spronk} that if $G$\/ is
compact, connected and non-abelian, then $(\Delta_G\times
\Delta_G)\Delta_{G\times G}$ is not a set of synthesis for
$A(G^4)$, where $G^4=G\times G\times G\times G$.

In this section, we show that for a general Lie group,
$(\Delta_G\times \Delta_G)\Delta_{G\times G}$ is always a set of
local weak synthesis. As in the case of the anti-diagonal, this is
done by showing first that $(\Delta_G\times
\Delta_G)\Delta_{G\times G}$ has the cone property.

\begin{lem}\label{L:product-submanifold}
Let $G$\/ be an $n$-dimensional Lie group. Then $(\Delta_G\times \Delta_G)\Delta_{G\times G}$ is a closed smooth submanifold of $G^4$ of dimension\/ $3n$, and is diffeomorphic to $G^3$.
\end{lem}

\begin{proof}
An element in $\Delta_G\times \Delta_G$ is of the form
$(s,s,v,v)$. An element in $\Delta_{G\times G}$ is of the form
$(w,z,w,z)$. Hence an element in $(\Delta_G\times
\Delta_G)\Delta_{G\times G}$ is of the form $(sw,sz,vw,vz)$.
Define
\begin{equation*}
M=\{(g_1,g_2,g_3,g_4)\in G^4 \mid g_1 g_3^{-1} g_4 g_2^{-1} =e\}.
\end{equation*}
Since $sw(vw)^{-1}vz(sz)^{-1}=e$, we have $(\Delta_G\times
\Delta_G)\Delta_{G\times G}\subset M$. Conversely, any element in
$M$\/ is of the form $(g_1,g_1g_3^{-1}g_4,g_3,g_4)$. If we set
$s=g_1$, $v=g_3$, $w=e$, and $z=g_3^{-1}g_4$, then we get
$(g_1,g_1g_3^{-1}g_4,g_3,g_4)=(sw,sz,vw,vz)$ and thus $M\subset
(\Delta_G\times \Delta_G)\Delta_{G\times G}$.  It follows that
$M=(\Delta_G\times \Delta_G)\Delta_{G\times G}$.

Now consider the smooth map $f: G^4\rightarrow G$\/ given by
\begin{equation}\label{eq: f}
f(g_1,g_2,g_3,g_4)=g_1 g_3^{-1} g_4 g_2^{-1}.
\end{equation}
Note that $M=f^{-1}(e)$ is closed in $G^4$.  Clearly, $f$\/ is
surjective. For fixed $g=(g_1,g_2,g_3,g_4)\in G^4$, consider the
right multiplication map $\rho:G\rightarrow G$\/ given by
\begin{equation*}
\rho(h)=hg_3^{-1} g_4 g_2^{-1}.
\end{equation*}
Since $\rho$\/ is a diffeomorphism, its derivative
$(D\rho)_{g_1}:T_{g_1}G\rightarrow T_{f(g)}G$\/ is bijective for
every $g_1\in G$. Since the image of $(D\rho)_{g_1}$ is contained
in the image of the derivative $(Df)_g: T_g G^4\rightarrow
T_{f(g)}G$, we conclude that $(Df)_g$\/ must be surjective.  It
follows that every element of $G$\/ is a regular value of $f$. In
particular, $e\in G$\/ is a regular value of $f$, and
consequently $M=f^{-1}(e)$ is a smooth submanifold of $G^4$ of
dimension equal to $\dim(G^4)-\dim G = 4n-n=3n$.

Finally, a diffeomorphism $\varphi : G^3\rightarrow (\Delta_G\times \Delta_G)\Delta_{G\times G}$ is given by
\begin{equation}\label{eq: varphi}
\varphi(g_1,g_2,g_3) = (g_1,g_1g_2^{-1}g_3,g_2,g_3).
\end{equation}
Note that $\varphi^{-1}$ is just the projection map
$(g_1,g_2,g_3,g_4)\mapsto (g_1,g_3,g_4)$ restricted to $(\Delta_G\times \Delta_G)\Delta_{G\times G}$.
\end{proof}

\begin{lem}\label{L:product-cone property}
Let $G$\/ be an $n$-dimensional Lie group. Then the closed
submanifold\/ $(\Delta_G\times \Delta_G)\Delta_{G\times G}$ has the cone
property in $G^4$.
\end{lem}

\begin{proof}
We proceed as in the proof of Lemma~\ref{L:cone property-anti diagonal}. Using the notation in Definition~\ref{defn: cone}, we set $H=G^4$, $E=M=(\Delta_G\times \Delta_G)\Delta_{G\times G}$, $m=3\dim G=3n$, and $U_x=G^4$
for every $x\in (\Delta_G\times \Delta_G)\Delta_{G\times G}$.  Let $\mathfrak{g}$ denote the Lie
algebra of $G$.  Let $W$\/ be a bounded open subset of $\mathfrak{g}\cong
\R^n$ containing $0$ such that the exponential map restricted to
$W$, $\exp :W\rightarrow \exp(W)\subset G$, is a diffeomorphism.
We set
$$W_x^0=W_x=W^3=W\times W\times W \subset \mathfrak{g}\times\mathfrak{g}\times\mathfrak{g} \cong \R^{3n},$$
and set $A_x ={\rm Aff}(G^4)$ for
every $x\in (\Delta_G\times \Delta_G)\Delta_{G\times G}$.  Given $p\in W$, we write $\tilde{p}=\exp(p)$ for short.
We define $\psi_x: W^3 \rightarrow {\rm Aff}(G^4)$,
also independent of $x$, as follows:  for $t = (p,q,r)\in W^3$ and
$(g_1,g_2,g_3,g_4)\in G^4$,
\begin{eqnarray*}
&&\psi_x(t)(g_1,g_2,g_3,g_4)
=(\tilde{p}g_1,\, \tilde{p} g_2 \tilde{q}\tilde{p}^{-1},\, \tilde{r}g_3,\, \tilde{r}g_4\tilde{q}\tilde{p}^{-1})\\
&&\hspace{20pt} =
(\tilde{p}g_1,\, (\tilde{q}\tilde{p} ^{-1})^{-1}(\tilde{q}g_2) \tilde{q}\tilde{p}^{-1},\, \tilde{r}g_3,\, (\tilde{q}\tilde{p} ^{-1})^{-1}(\tilde{q}\tilde{p} ^{-1}\tilde{r}g_4)\tilde{q}\tilde{p}^{-1}).
\end{eqnarray*}
Each $\psi_x(t)$ is the composition of a left multiplication by $(\tilde{p},\tilde{q},\tilde{r},\tilde{q}\tilde{p}^{-1}\tilde{r})$ and an inner automorphism $({\rm id}_G \times {\rm Inn}(\tilde{q}\tilde{p}^{-1}) \times {\rm id}_G \times {\rm Inn}(\tilde{q}\tilde{p}^{-1})) \in {\rm Aut}(G^4)$, and hence is an affine transformation of $G^4$.
Clearly, $\psi_x$ is smooth and $\psi_x(0)={\rm id}_{G^4}$.  For every $y= (g_1,g_1g_3^{-1}g_4,g_3,g_4)\in (\Delta_G\times \Delta_G)\Delta_{G\times G}$ and $t=(p,q,r)\in W^3$, we have
\begin{eqnarray*}
f(\psi_x(t)y)&=&f(\tilde{p}g_1,\, \tilde{p} g_1g_3^{-1}g_4 \tilde{q}\tilde{p}^{-1},\, \tilde{r}g_3,\, \tilde{r}g_4\tilde{q}\tilde{p}^{-1})\\
&=& (\tilde{p}g_1) (\tilde{r}g_3)^{-1} ( \tilde{r}g_4\tilde{q}\tilde{p}^{-1}) (\tilde{p} g_1g_3^{-1}g_4 \tilde{q}\tilde{p}^{-1})^{-1}=e,
\end{eqnarray*}
where $f:G^4\rightarrow G$\/ is the smooth map defined in (\ref{eq: f}).  It follows that
\begin{equation*}
\psi_x(t)y \,\in f^{-1}(e)=(\Delta_G\times \Delta_G)\Delta_{G\times G}.
\end{equation*}

It only remains to prove that $W^3\rightarrow \psi_x(W^3)y$\/ is a diffeomorphism for fixed $y=(g_1,g_1g_3^{-1}g_4,g_3,g_4)\in (\Delta_G\times \Delta_G)\Delta_{G\times G}$.  Let $\varphi : G^3 \rightarrow (\Delta_G\times \Delta_G)\Delta_{G\times G}$ be the diffeomorphism defined in (\ref{eq: varphi}).  It is enough to show that the map $\xi : (\exp (W))^3 \rightarrow \varphi^{-1}(\psi_x(W^3)y)$\/ given by
\begin{equation*}
\xi (\tilde{p},\tilde{q},\tilde{r}) =  (\tilde{p}g_1,\, \tilde{r}g_3,\, \tilde{r}g_4\tilde{q}\tilde{p}^{-1})
\end{equation*}
is a diffeomorphism.  Clearly, $\xi$\/ extends to a smooth map $\hat{\xi} : G^3\rightarrow G^3$.  It is easy to check that $\hat{\xi}$\/
is bijective and the inverse map is given by
\begin{equation*}
\hat{\xi}^{-1}(z_1,z_2,z_3) = (z_1g_1^{-1},\, g_4^{-1}g_3 z_2^{-1}z_3 z_1g_1^{-1},\, z_2g_3^{-1}),
\end{equation*}
which is smooth as well.  Hence $\hat{\xi}$\/ is a diffeomorphism, and so is its restriction $\xi=\hat{\xi}|_{(\exp(W))^3}$.
\end{proof}

\begin{thm}\label{T:smooth weak syn-product}
Let $G$ be an $n$-dimensional Lie group. Then:\\
{\rm (i)} $(\Delta_G\times \Delta_G)\Delta_{G\times G}$ is a set of local smooth synthesis for $A(G^4)$;\\
{\rm (ii)} $(\Delta_G\times \Delta_G)\Delta_{G\times G}$ is a set of local weak synthesis for
$A(G^4)$ of degree at most\/ $[\frac{3n}{2}]+1$.
\end{thm}

\begin{proof}
These follow immediately from Theorem \ref{T:cone prop-smooth weak syn} and Lemma \ref{L:product-cone property}.
\end{proof}

\section{Projection theorem for sets of smooth and weak
synthesis}\label{S:Proj-smooth-weak}


In this section, we prove the projection theorem for sets of weak
and smooth synthesis.
Let $G$\/ be a locally compact group, and let $\fA(G)$ be a regular
Banach algebra of continuous functions on $G$\/ which is closed
under right translations and such that for any $f\in \fA(G)$ we
have
\begin{align*}
&\bullet\; \| R_sf\|_\fA=\|f\|_\fA \text{ for any }s \in
G,\\
&\bullet\; s\mapsto R_sf\text{ is continuous}.
\end{align*}

If $K$\/ is a compact subgroup of $G$\/ we let
\[
\fA(G:K)=\{f\in\fA(G)\mid
R_k f=f\hspace{1pt}\text{ for each }k\in K\},
\]
which is a closed subalgebra of $\fA(G)$ whose elements
are constant on left cosets of $K$.  We let $G/K$
denote the space of left cosets with the quotient topology.
We define two maps
\begin{align*}
P: \fA(G)\to \fA(G), \quad  Pf=\int_K R_k f \hspace{1pt}dk, \\
M: \fA(G:K)\to C_b(G/K),\quad  Mf(sK)=f(s),
\end{align*}
where $C_b(X)$ denotes the space of bounded continuous functions on a
topological space $X$. The map $P$\/ is to be regarded as a Bochner
integral over the normalized Haar measure on $K$; its range is
$\fA(G:K)$ and $P$\/ is a (completely) contractive projection. The
map $M$\/ is well-defined by comments above, and its range consists
of continuous functions since $\fA(G:K)\subset C_b(G:K)$. We note
that $M$\/ is an injective homomorphism and denote its range by
$\fA(G/K)$.  We assign a norm to $\fA(G/K)$ in such a way that $M$\/
is an isometry. We finally define two maps
\[
N=M^{-1}:\fA(G/K)\to \fA(G), \ \ \Gamma=M\circ P:\fA(G)\to\
\fA(G/K)
\]
so that $N$\/ is an isometric homomorphism and $\Gamma$ is a quotient
map.

The following theorem demonstrates the relationship between (local)
weak synthesis for $\fA(G)$ and $\fA(G/K)$. Its proof follows directly
from \cite[Theorem 1.4]{forrest-samei-spronk} (see also \cite[Corollary 1.5]{forrest-samei-spronk}).

\begin{thm}\label{T:coset-weak syn}
Let $G$ be a compact group, and let $\fT(G)$ denote the space of trigonometric polynomials on $G$.
Let $\fA(G)$ be as above and
additionally satisfy that $\fT(G)\fA(G) \subseteq \fA(G)$. If $E$
is a closed subset of $G/K$, let $$E^*=\{s\in G \mid sK\in E\}.$$ Then,
for every $d\in \N$, $E$ is a set of\/ $($local\/$)$ weak synthesis for
$\fA(G/K)$ of degree at most $d$ if and only if $E^*$ is a set of\/
$($local\/$)$ weak synthesis for $\fA(G)$ of degree at most $d$ .


\end{thm}

Now let $H$\/ be a subgroup of a Lie group $G$.  Define
$$C^{\infty}(G:H)=\{f\in C^{\infty}(G) \mid f(sh)=f(s),\ \forall
s\in G,\ \forall h\in H\}.$$

\begin{lem}\label{L:smooth-subgroup}
Let $G$ be a Lie group, and let $H$ be a compact
subgroup of $G$.  If $f\in C^{\infty}(G)$, then for $s\in G$ the
Haar integral
\begin{equation*}
(P_H f)(s) = \int_H f(sr) dr
\end{equation*}
defines an element in $C^{\infty}(G:H)$.
\end{lem}

\begin{proof}
It is well-known that $H$\/ is also a Lie group. Let $n=\dim G$\/ and $k=\dim H$.  Since $H$\/ is compact, there exists a smooth orientation $k$-form $\omega$\/ on $H$\/ that is invariant under both left and right multiplications.  By the definition of Haar integral, we have
\[
\int_H \omega = 1 \quad \textrm{and} \quad \int_H f(sr) dr = \int_H f(sr) \omega .
\]
Since $\omega$\/ is invariant under left multiplication, it is clear that $(P_H f)(sh)=(P_H f)(s)$ for all pairs $s\in G$ and $h\in H$.
It remains to prove that $P_H f \in C^{\infty}(G)$.

Let $U$\/ be an open chart of $G$.  Choose a coordinate system
$s=(x^1,\dots,x^n)$ on $U$.
Let $\{V_i \mid i=1,\dots, m\}$ be oriented open charts of $H$\/ that cover the compact manifold $H$, and let $\{\psi_i : H \rightarrow [0,1]  \mid i=1,\dots, m\}$ be a smooth partition of unity subordinate to this cover.  Choose a coordinate system $y_i=(y_i^1,\dots,y_i^k)$ on each open subset $V_i$ of $H$\/ such that
\[
\omega|_{V_i} =  \lambda_i(y_i)\, dy_i^1 \wedge \dots \wedge dy_i^k,
\]
where $\lambda_i(y_i)=\lambda_i(y_i^1,\dots,y_i^k)$ is a positive smooth function of $y_i\in V_i$.

Let $F:G\times H\rightarrow \R$ be the smooth function given by $F(s,r)=f(\mu(s,r))=f(sr)$, where $\mu:G\times G \rightarrow G$\/ is the group multiplication map.  Now we can write
\[
F|_{U\times V_i}(s,y_i) = \tilde{F}_i (x^1,\dots,x^n,y_i^1,\dots,y_i^k),
\]
where $\tilde{F}_i$ is a smooth function defined on the open subset of $\R^n\times \R^k$ that corresponds to $U\times V_i$ under the above coordinate systems.  It follows that
\begin{equation*}
(P_H f)|_{U} (s) =
\sum_{i=1}^m \int_{V_i} \psi_i(y_i) \lambda_i(y_i) \tilde{F}_i (x^1,\dots,x^n,y_i^1,\dots,y_i^k) \, dy_i^1  \dots dy_i^k .
\end{equation*}
Each multiple Riemann integral in the above sum is a smooth
function of the multivariable $s=(x^1,\dots,x^n)$, and this shows
that $P_H f \in C^{\infty}(G)$.
\end{proof}

\begin{thm}\label{T:coset-smooth syn}
Let $G$ be a compact Lie group, let $\fA(G)$ be as in Theorem\/
$\ref{T:coset-weak syn}$, and additionally contain $\D(G)$.  If
$E$ is a closed subset of $G/K$, let $E^*=\{s\in G \mid sK\in E\}$.
Then $E$ is a set of smooth synthesis for $\fA(G/K)$ if and only
if $E^*$ is a set of smooth synthesis for $\fA(G)$.
\end{thm}

\begin{proof}
This follows from \cite[Theorem 1.4]{forrest-samei-spronk} and Lemma \ref{L:smooth-subgroup}.
\end{proof}

\section{Alternative proofs in the case of a compact Lie group}
\label{S:Alternative}

In this section, we apply Theorem \ref{T:coset-smooth syn} to
present alternative proofs of the smooth and weak synthesis results
of Section \ref{Diag type-smooth-weak} for compact Lie groups
(Theorems \ref{T:smooth weak syn-anti diagonal} and \ref{T:smooth
weak syn-product}).

Let $G$\/ be a {\it compact}\/ Lie group for the remainder of this
section. We use $G\times G$ in place of $G$, and
$K=\Delta=\{(s,s)\mid s\in G\}$ in the setup of Section \ref{S:Proj-smooth-weak}.
Since the map
\begin{equation}\label{eq:ggmoddel}
(G\times G)/\Del \to G, \quad (s,e)\Delta \mapsto s
\end{equation}
is a homeomorphism, we identify the coset space with $G$. We
observe that in this case the map $P:A(G\times G) \to
A(G\times G)$ satisfies
\[
Pw(s,t)=\int_G w(sr,tr)dr=\int_G w(st^{-1}r,r)dr
\]
and the map $M:A(G\times G:\Del)\to C(G)$ satisfies
\begin{equation*}
Mw(s)=w(s,e).
\end{equation*}
The map $\Gamma=M\circ P$, from Section~\ref{S:Proj-smooth-weak}, can
be regarded as a `twisted' convolution, for if $A(G\cross G)$
contains an elementary function $f\cross g$, then for $s\in G$
\[
\Gamma(f\cross g)(s)=\int_G (f\cross g)(st,t)dt
=\int_G f(st)g(t)dt=f*\check{g}(s).
\]
We denote the image of $\Gamma$ by $A_\Del(G)$. We endow
$A_\Delta(G)$ with the norm which makes $\Gamma$ a quotient map.
We also note that
\[
N:A_\Delta(G) \rightarrow A(G\cross G:\Delta), \quad
Nu(s,t)=u(st^{-1})
\]
is an isometry. As in \cite[Theorem 2.6]{forrest-samei-spronk}, if
we repeat the procedure above we obtain
$$A_{\Del^2}(G)=\Gamma(A_\Del(G\times G)).$$


We can do a similar construction with the anti-diagonal
$\check{\Del}=\{(s,s^{-1})\mid s\in G\}$. We let $G\cross
G/\check{\Del}$ denote the set of equivalence classes modulo the
equivalence relation $(s',t')\sim (s,t)$ if and only if
$(s^{-1}s' , t' t^{-1})\in\check{\Del}$, so that $G\cross G/\check{\Del}$,
with the quotient topology, is homeomorphic to $G$\/ via
$(s,t)\mapsto st$.  We let
\[
A(G\cross G:\check{\Del})=\{u\in A(G\times G) \mid r\diamond u=u \
\text{for all} \ r\in G\},
\]
where $(r\diamond u)(s,t)=u(sr,r^{-1}t)$. Similarly as above,
$A(G\cross G:\check{\Del})$ is a closed subalgebra of $A(G\times
G)$. Also the map
\[
\check{\Gamma}:A(G\times G) \rightarrow C(G) , \quad
\check{\Gamma} w(s)=\int_G w(st,t^{-1})dt
\]
is surjective, and is injective on $A(G\cross G:\check{\Del})$.
We denote the image of $\check{\Gamma}$ by $A_\gamma(G)$. We endow
$A_\gamma(G)$ with the norm which makes $\check{\Gamma}$ a
quotient map.  We also note that
\[
\check{N}:A_\gamma(G) \rightarrow A(G\cross G:\check{\Del}) ,
\quad \check{N}u(s,t)=u(st)
\]
is an isometry.

\begin{lem}\label{lem: GxG}
Let $G$ be a compact Lie group.  If $F\in C^{\infty}(G\times G)$, then for $s\in G$ the Haar integral
\[
\hat{F}(s) = \int_G F(s,r)dr
\]
defines an element in $C^{\infty}(G)$.
\end{lem}

\begin{proof}
We argue similarly as in the proof of Lemma~\ref{L:smooth-subgroup}.  Let $n=\dim G$, and let $\omega$\/ be a smooth bi-invariant orientation $n$-form on $G$\/ such that
\[
\int_G \omega = 1 \quad \textrm{and} \quad \int_G F(s,r) dr = \int_G F(s,r) \omega.
\]
Let $\{U_i \mid i=1,\dots, m\}$ be oriented open charts of $G$\/ that cover the compact manifold $G$, and let $\{\psi_i : G \rightarrow [0,1]  \mid i=1,\dots, m\}$ be a smooth partition of unity subordinate to this cover.  Choose a coordinate system $x_i=(x_i^1,\dots,x_i^n)$ on each open subset $U_i$ of $G$\/ such that
\[
\omega|_{U_i} =  \lambda_i(x_i)\, dx_i^1 \wedge \dots \wedge dx_i^n,
\]
where $\lambda_i(x_i)=\lambda_i(x_i^1,\dots,x_i^n)$ is a positive smooth function of $x_i\in U_i$.

Given $s\in G$, we have $s=(x_j^1,\dots,x_j^n)\in U_j$ for some $j \in\{1,\dots,m\}$.  We fix one such $j$\/ and write
\[
F|_{U_j\times U_i}(s,x_i) = \tilde{F}_{j,i} (x_j^1,\dots,x_j^n,x_i^1,\dots,x_i^n),
\]
where $\tilde{F}_{j,i}$ ($i=1,\dots,m$) is a smooth function defined on the open subset of $\R^n\times \R^n$ that corresponds to $U_j\times U_i\subset G\times G$\/ under the above coordinate systems.  It follows that
\begin{equation*}
\hat{F}|_{U_j} (s) =
\sum_{i=1}^m \int_{U_i} \psi_i(x_i) \lambda_i(x_i) \tilde{F}_{j,i} (x_j^1,\dots,x_j^n,x_i^1,\dots,x_i^n) \, dx_i^1  \dots dx_i^n .
\end{equation*}
Each multiple Riemann integral in the above sum is a smooth function of the multivariable $s=(x_j^1,\dots,x_j^n)$, and this shows that $\hat{F} \in C^{\infty}(G)$.
\end{proof}

\begin{lem}\label{L:smooth-diagonal}
Let $G$ be a compact Lie group. Then:\\
{\rm (i)} $\Gamma\circ N = {\rm id}_{C(G)}$ and\/ $\check{\Gamma}\circ \check{N} = {\rm id}_{C(G)}$.\\
{\rm (ii)} If $f\in C^{\infty}(G)$, then\/ $Nf \in C^{\infty}(G\times G)$ and\/ $\check{N}f\in C^{\infty}(G\times G)$.\\
{\rm (iii)}  If $f\in C^{\infty}(G\times G)$, then\/ $\Gamma f \in C^{\infty}(G)$ and\/ $\check{\Gamma}f\in C^{\infty}(G)$.
\end{lem}

\begin{proof}
Part (i) is obvious.  For part (ii), let $\mu:G\times G\rightarrow G$, $(s,t)\mapsto st$, denote the multiplication map, and let ${\rm inv}_G:G\rightarrow G$, $g\mapsto g^{-1}$, denote the inverse map.  Since $\mu$\/ and ${\rm inv}_G$ are both smooth, so are the compositions
$\check{N}f=f\circ \mu$\/ and $Nf=f\circ\mu\circ({\rm id}_G \times {\rm inv}_G)$.
For part (iii), we apply Lemma~\ref{lem: GxG} to the smooth map $F(s,r)=(f\circ(\mu\times {\rm id}_G))(s,r)=f(sr,r)$ to conclude that $\hat{F}=\Gamma f$\/ is smooth.  Similarly, we apply Lemma~\ref{lem: GxG} to $F(s,r)=(f\circ(\mu\times{\rm inv}_G))(s,r)=f(sr,r^{-1})$ to conclude that $\hat{F}=\check{\Gamma}f$\/ is smooth.
\end{proof}

We are now ready to state the main result of this section.

\begin{thm}\label{T:smooth weak syn-anti diagonal-product-alternative}
Let $G$ be an $n$-dimensional compact Lie group. Then:\\
{\rm (i)}  $\check{\Delta}_G$ is a set of smooth synthesis for $A(G\times G)$.\\
{\rm (ii)} $\check{\Delta}_G$ is a set of weak synthesis for $A(G\times G)$ of degree at most\/ $[\frac{n}{2}]+1$.\\
{\rm (iii)} $(\Delta_G\times \Delta_G)\Delta_{G\times G}$ is a set of smooth synthesis for $A(G^4)$.\\
{\rm (iv)} $(\Delta_G\times \Delta_G)\Delta_{G\times G}$ is a set of weak synthesis for
$A(G^4)$ of degree at most\/ $[\frac{3n}{2}]+1$.
\end{thm}

\begin{proof}
(i) It is noted in the proof of Lemma 2.3 in \cite{FSS2} that
$$I(\check{\Delta}_G)=\overline{\spn} \{u\check{N}f \mid u\in A(G\times G) \ \text{and}\
f\in I_{A_\gamma(G)}(\{e\}) \}.$$ Let us first show that $\{e\}$
is a set of smooth synthesis for $A_\gamma(G)$, i.e.
\begin{equation}\label{eq:J-1}
J^\D_{A_\gamma(G)}(\{e\})=I_{A_\gamma(G)}(\{e\}).
\end{equation}
Since $D(G)\cap A(G)$ is dense in $A(G)$, for every $u\in
I_{A_\gamma(G)}(\{e\})$ there is a sequence $\{ u_n \} \subset
D(G)\cap A(G)$ which converges to $u$. We note that
\[
|u_n(e)|=|u_n(e)-u(e)|\leq \|u_n-u\|_\infty \leq \|u_n-u\|_{A_\gamma(G)}
\overset{n\to\infty}{\longrightarrow}0.
\]
Thus if $u_n'=u_n-u_n(e)1$, then $\{ u_n'\}\subset
J^\D_{A_\gamma(G)}(\{e\})$ with $\lim_{n\to\infty}
\|u_n'-u\|_{A_\gamma(G)}=0$. Hence (\ref{eq:J-1}) holds. Now it follows from
Lemma \ref{L:smooth-diagonal}(ii)
that $\check{\Del}_G$ is a set of smooth synthesis for $A(G\times G)$.\\
(ii) This follows from part (i) and Theorem \ref{T:smooth-weak synthesis}.\\
(iii) It is shown in \cite[Theorems 1.4 and 2.6]{forrest-samei-spronk} that
$$I_{A_\Del(G\times G)}(\Del_G)=\overline{\spn} \{uNf \mid u\in A_\Del(G\times G) \ \text{and}\
f\in I_{A_{\Del^2}(G)}(\{e\}) \}.$$ We can see, similar to (\ref{eq:J-1})
that $\{e\}$ is a set of smooth synthesis for $A_{\Del^2}(G)$,
i.e.
\begin{equation}\label{eq:J-2}
J^\D_{A_{\Del^2}(G)}(\{e\})=I_{A_{\Del^2}(G)}(\{e\}).
\end{equation}
Thus, from (\ref{eq:J-2}) and Lemma \ref{L:smooth-diagonal}(ii), it follows
that $\Del_G$ is a set of smooth synthesis for $A_\Del(G\times
G)$, and so, $(\Delta_G \times \Del_G)\Del_{G\times G}$  is a set
of smooth synthesis for $A(G\times G \times G \times G)$ from
Theorem \ref{T:coset-smooth syn} and
\cite[Theorems 1.4 and 2.6]{forrest-samei-spronk}.\\
(iv) This follows from part (iii) and Theorem \ref{T:smooth-weak synthesis}.
\end{proof}


\subsection*{Acknowledgments}  The authors thank Dragomir \v{Z}.
Dokovi\'c for very helpful discussions. The first author was
partially supported by CFI, NSERC and OIT grants and the second
author was partially supported by an NSERC Postdoctoral Fellowship.
The second author also would like to thank Nico Spronk for reading
the earlier draft of this article and giving valuable comments.


\begin{thebibliography}{99}

\bibitem{D}  H. G. Dales, \textit{Banach algebras and automatic
continuity}, Oxford University Press, New York, 2000.

\bibitem{domar}  Y. Domar,
\textit{On the spectral synthesis problem for\/ $(n-1)$-dimensional subsets of\/
$\R^n,$ $n\geq 2$},
Ark. Mat. \textbf{9} (1971), 23--37.

\bibitem{D2}  Y. Domar, \textit{A $C^\infty$ curve of spectral
non-synthesis}, Mathematika \textbf{24} (1977), 189--192.

\bibitem{D3}  Y. Domar, \textit{On spectral synthesis in\/ $\R^n$, $n\geq 2$},
Lecture Notes in Math., vol. 779, Springer, Berlin, 1980, pp. 46--72.

\bibitem {Em}  P. Eymard, \textit{L'alg\`{e}bre de Fourier d'un groupe
localement compact}, Bull. Soc. Math. France \textbf{92} (1964),
181--236.

\bibitem{FR}  B. E. Forrest and V. Runde, \textit{Amenability and weak amenability of
the Fourier algebra}, Math. Z. \textbf{250} (2005), 731--744.

\bibitem{forrest-samei-spronk}  B. E. Forrest, E. Samei and N. Spronk,
\textit{Convolutions on compact groups and Fourier algebras of coset spaces},
arXiv:0705.4277.

\bibitem{FSS2}  B. E. Forrest, E. Samei and N. Spronk,
\textit{Weak amenability of Fourier algebras on compact groups},
arXiv:0808.1858.

\bibitem{guo}  K. Guo,
\textit{A remark on the spectral synthesis property for hypersurfaces of\/ $\R^n$},
Proc. Amer. Math. Soc. \textbf{121} (1994), 185--192.

\bibitem{helgason}  S. Helgason,
\textit{Differential geometry, Lie groups, and symmetric spaces},
Pure Appl. Math., vol. 80, Academic Press, Inc., New York, 1978.

\bibitem{H1}  C. Herz, \textit{Harmonic synthesis for subgroups}, Ann.
Inst. Fourier (Grenoble) \textbf{23} (1973), 91--123.

\bibitem{Kan}  E. Kaniuth, \textit{Weak spectral synthesis in commutative Banach algebras},
J. Funct. Anal. \textbf{254} (2008), 987--1002.

\bibitem{KM}  W. Kirsch and D. M\"uller,
\textit{On the synthesis problem for orbits of Lie groups
in\/ $\R^n$}, Ark. Mat. \textbf{18} (1980), 145--155.

\bibitem{ludwig-turowska}  J. Ludwig and L. Turowska,
\textit{Growth and smooth spectral synthesis in the Fourier
algebras of Lie groups}, Studia Math. \textbf{176} (2006),
139--158.

\bibitem{Meaney}  C. Meaney, \textit{On the failure of spectral synthesis for compact semisimple Lie groups}, J. Funct. Anal. \textbf{48} (1982), 43--57.


\bibitem{muller}  D. M\"uller,
\textit{On the spectral synthesis problem for hypersurfaces of\/ $\R^n$},
J. Funct. Anal. \textbf{47} (1982), 247--280.

\bibitem{Ric}  C. E. Rickart, \textit{General theory of Banach algebras},
Van Nostrand, Princeton, NJ, 1960.

\bibitem{TT}  M. Takesaki and N. Tatsuuma, \textit{Duality and subgroups II},
J. Funct. Anal. \textbf{11} (1972), 184--190.

\bibitem{V}  N. Th. Varopoulos, \textit{Spectral synthesis on
spheres}, Proc. Cambridge Philos. Soc. \textbf{62} (1966), 379--387.

\bibitem{W}  C. R. Warner, \textit{Weak spectral synthesis}, Proc.
Amer. Math. Soc. \textbf{99} (1987), 244--248.

\end{thebibliography}
\end{document}